# ON THE $p$−LOGARITHMIC AND $\alpha$−POWER DIVERGENCE MEASURES IN INFORMATION THEORY

S.S. DRAGOMIR

ABSTRACT. In this paper we introduce the concepts of $p$−logarithmic and $\alpha$−power divergence measures and point out a number of basic results.

## 1. INTRODUCTION

One of the important issues in many applications of Probability Theory is finding an appropriate measure of *distance* (or *difference* or *discrimination* ) between two probability distributions. A number of divergence measures for this purpose have been proposed and extensively studied by Jeffreys [2], Kullback and Leibler [3], Rényi [4], Havrda and Charvat [5], Kapur [6], Sharma and Mittal [7], Burbea and Rao [8], Rao [9], Lin [10], Csiszár [11], Ali and Silvey [13], Vajda [14], Shioya and Da-te [48] and others (see for example [6] and the references therein).

These measures have been applied in a variety of fields such as: anthropology [9], genetics [15], finance, economics, and political science [16], [17], [18], biology [19], the analysis of contingency tables [20], approximation of probability distributions [21], [22], signal processing [23], [24] and pattern recognition [25], [26].

Assume that a set $\chi$ and the $\sigma$−finite measure $\mu$ are given. Consider the set of all probability densities on $\mu$ to be $\Omega := \left\{ p | p : \chi \to \mathbb{R}, \ p(x) \geq 0, \ \int_\chi p(x) \, d\mu(x) = 1 \right\}$. The Kullback-Leibler divergence [3] is well known among the information divergences. It is defined as:

$$(1.1) \qquad D_{KL}(p, q) := \int_\chi p(x) \log \left[ \frac{p(x)}{q(x)} \right] d\mu(x), \ \ p, q \in \Omega,$$

where log is to base 2.

In Information Theory and Statistics, various divergences are applied in addition to the Kullback-Leibler divergence. These are the: *variation distance* $D_v$, *Hellinger distance* $D_H$ [41], $\chi^2$−*divergence* $D_{\chi^2}$, $\alpha$−*divergence* $D_\alpha$, *Bhattacharyya distance* $D_B$ [42], *Harmonic distance* $D_{Ha}$, *Jeffrey's distance* $D_J$ [2], *triangular discrimination* $D_\Delta$ [36], etc... They are defined as follows:

$$(1.2) \qquad D_v(p, q) := \int_\chi |p(x) - q(x)| \, d\mu(x), \ \ p, q \in \Omega;$$

$$(1.3) \qquad D_H(p, q) := \int_\chi \left| \sqrt{p(x)} - \sqrt{q(x)} \right| d\mu(x), \ \ p, q \in \Omega;$$

$$(1.4) \quad D_{\chi^2}(p, q) := \int_\chi p(x) \left[\left(\frac{q(x)}{p(x)}\right)^2 - 1\right] d\mu(x), \ p, q \in \Omega;$$

$$(1.5) \quad D_\alpha(p, q) := \frac{4}{1-\alpha^2}\left[1 - \int_\chi [p(x)]^{\frac{1-\alpha}{2}} [q(x)]^{\frac{1+\alpha}{2}} d\mu(x)\right], \ p, q \in \Omega;$$

$$(1.6) \quad D_B(p, q) := \int_\chi \sqrt{p(x) q(x)} d\mu(x), \ p, q \in \Omega;$$

$$(1.7) \quad D_{Ha}(p, q) := \int_\chi \frac{2 p(x) q(x)}{p(x) + q(x)} d\mu(x), \ p, q \in \Omega;$$

$$(1.8) \quad D_J(p, q) := \int_\chi [p(x) - q(x)] \ln\left[\frac{p(x)}{q(x)}\right] d\mu(x), \ p, q \in \Omega;$$

$$(1.9) \quad D_\Delta(p, q) := \int_\chi \frac{[p(x) - q(x)]^2}{p(x) + q(x)} d\mu(x), \ p, q \in \Omega.$$

For other divergence measures, see the paper [6] by Kapur or the book on line [43] by Taneja. For a comprehensive collection of preprints available on line, see the RGMIA web site http://rgmia.vu.edu.au/papersinfth.html

Csiszár $f-$divergence is defined as follows [11]

$$(1.10) \quad D_f(p, q) := \int_\chi p(x) f\left[\frac{q(x)}{p(x)}\right] d\mu(x), \ p, q \in \Omega,$$

where $f$ is convex on $(0, \infty)$. It is assumed that $f(u)$ is zero and strictly convex at $u = 1$. By appropriately defining this convex function, various divergences are derived. All the above distances $(1.1) - (1.9)$, are particular instances of Csiszár $f-$divergence. There are also many others which are not in this class (see for example [6] or [43]). For the basic properties of Csiszár $f-$divergence see [44]-[46].

In [47], Lin and Wong (see also [10]) introduced the following divergence

$$(1.11) \quad D_{LW}(p, q) := \int_\chi p(x) \log\left[\frac{p(x)}{\frac{1}{2}p(x) + \frac{1}{2}q(x)}\right] d\mu(x), \ p, q \in \Omega.$$

This can be represented as follows, using the Kullback-Leibler divergence:

$$D_{LW}(p, q) = D_{KL}\left(p, \frac{1}{2}p + \frac{1}{2}q\right).$$

Lin and Wong have established the following inequalities

$$(1.12) \quad D_{LW}(p, q) \leq \frac{1}{2} D_{KL}(p, q);$$

$$(1.13) \quad D_{LW}(p, q) + D_{LW}(q, p) \leq D_v(p, q) \leq 2;$$

$$(1.14) \quad D_{LW}(p, q) \leq 1.$$

In [48], Shioya and Da-te improved $(1.12) - (1.14)$ by showing that

$$D_{LW}(p,q) \leq \frac{1}{2}D_v(p,q) \leq 1.$$

For classical and new results in comparing different kinds of divergence measures, see the papers [2]-[48] where further references are given.

## 2. Some New Divergence Measures

We define the $p-$Logarithmic means by (see [1, p. 346])

$$
(2.1) \quad L_p(a,b) = \begin{cases} \left[\frac{b^{p+1}-a^{p+1}}{(p+1)(b-a)}\right]^{\frac{1}{p}}, & p \neq -1, 0, \\ \frac{b-a}{\ln b - \ln a}, & p = -1, \quad a \neq b,\ a,b > 0 \\ \frac{1}{e}\left(\frac{b^b}{a^a}\right)^{\frac{1}{b-a}}, & p = 0, \end{cases}
$$

$$L_p(a,a) = a.$$

Where convenient, $L_{-1}(a,b)-$the logarithmic mean, will be written as just $L(a,b)$. The case $p = 0$ is also called the *identric mean*, i.e., $L_0(a,b)$ and will be denoted by $I(a,b)$. Of course, we will also define $L_\infty(a,b) = \max\{a,b\}$ and $L_{-\infty} = \min\{a,b\}$ to complete the scale.

It is easily checked that the definitions in the above scale are consistent in the sense that $\lim_{p \to 0} L_p(a,b) = I(a,b)$ and $\lim_{p \to \pm\infty} L_p(a,b) = L_{\pm\infty}(a,b)$.

We define the $p-$logarithmic divergence measure, or simply the $L_p-$*divergence measure,* by

$$
(2.2) \quad D_{L_p}(q,r) = \begin{cases} \int_\chi \left[\frac{[q(x)]^{p+1}-[r(x)]^{p+1}}{(p+1)(q(x)-r(x))}\right]^{\frac{1}{p}} d\mu(x), & \text{if } p \neq -1, 0, \\ \int_\chi \left[\frac{q(x)-r(x)}{\ln q(x) - \ln r(x)}\right] d\mu(x), & \text{if } p = -1, \quad q, r \in \Omega \\ \frac{1}{e}\int_\chi \left[\frac{[q(x)]^{q(x)}}{[r(x)]^{r(x)}}\right]^{\frac{1}{q(x)-r(x)}} d\mu(x), & \text{if } p = 0, \end{cases}
$$

$$(2.3) \quad D_{+\infty}(q,r) := \int_\chi \max\{q(x), r(x)\} d\mu(x), \quad q, r \in \Omega,$$

$$(2.4) \quad D_{-\infty}(q,r) := \int_\chi \min\{q(x), r(x)\} d\mu(x), \quad q, r \in \Omega.$$

We observe that

$$(2.5) \quad D_{+\infty}(q,r) = \int_\chi \frac{q(x)+r(x)+|q(x)-r(x)|}{2} d\mu(x) = 1 + \frac{1}{2}D_v(q,r)$$

and similarly,

$$(2.6) \quad D_{-\infty}(q,r) = 1 - \frac{1}{2}D_v(q,r).$$

Since $L_p(a,b) = L_p(b,a)$ for all $a, b > 0$ and $p \in [-\infty, \infty]$, we can conclude that the $L_p-$*divergence measures are symmetrical.*

Now, if we consider the continuous mappings (which are not necessarily convex)

$$(2.7) \quad f_p(x) := \begin{cases} \left[\frac{x^{p+1}-1}{(p+1)(x-1)}\right]^{\frac{1}{p}}, & x \in (0,1) \cup (1,\infty),\ p \neq -1, 0; \\ \frac{x-1}{\ln x-1}, & x \in (0,1) \cup (1,\infty),\ p = -1; \\ \frac{1}{e} x^{\frac{x}{x-1}}, & x \in (0,1) \cup (1,\infty),\ p = 0; \\ 1 & \text{if } x = 1, \end{cases}$$

$$f_{+\infty}(x) = 1 + \frac{1}{2}|x-1|,$$
$$f_{-\infty}(x) = 1 - \frac{1}{2}|x-1|,$$

and taking into account that $L_p(a,b) = aL_p\left(1, \frac{b}{a}\right)$ for all $a, b > 0$ and $p \in [-\infty, \infty]$, we deduce that

$$(2.8) \quad D_{f_p}(q,r) = \int_\mathcal{X} q(x) f_p\left(\frac{r(x)}{q(x)}\right) d\mu(x)$$
$$= \int_\mathcal{X} q(x) L_p\left(\frac{r(x)}{q(x)}, 1\right) d\mu(x)$$
$$= \int_\mathcal{X} L_p(r(x), q(x)) d\mu(x) = D_{L_p}(q,r)$$

for all $q, r \in \Omega$, which shows that the $L_p-$divergence measure can be interpreted as $f-$Csiszár divergences for $f = f_p$, which are not necessarily convex.

The following result is well known in the theory of $p-$logarithmic means (see for example [1, p. 347]).

**Lemma 1.** *We have*

$$(2.9) \qquad\qquad L_{-2}(a,b) = G(a,b),$$
$$(2.10) \qquad\qquad L_{-\frac{1}{2}}(a,b) = \frac{1}{2}(A(a,b) + G(a,b)),$$
$$(2.11) \qquad\qquad L_1(a,b) = A(a,b),$$
$$(2.12) \qquad\qquad L_{-3}(a,b) = \left[H(a,b)G^2(a,b)\right]^{\frac{1}{3}},$$

*and the monotonicity property*

$$(2.13) \qquad L_{-\infty}(a,b) \leq L_r(a,b) \leq L_s(a,b) \leq L_{+\infty}(a,b)$$

*with equality iff $a = b$, where $-\infty \leq r < b \leq \infty$, and $A(a,b)$ is the arithmetic mean, $G(a,b)$ is the geometric mean and $H(a,b)$ is the harmonic mean of $a, b$. In particular, we have*

$$(2.14) \quad G(a,b) \leq L(a,b) \leq \frac{1}{2}[A(a,b) + G(a,b)] \leq I(a,b) \leq A(a,b)$$

*with equality iff $a = b$.*

Now, using (2.9) - (2.12), we observe that

$$(2.15) \quad D_{L_{-2}}(q,r) = \int_\mathcal{X} \sqrt{r(x)q(x)} d\mu(x)$$
$$= D_B(p,q) \quad \text{(Bhattacharyya distance)}$$

$$(2.16) \quad D_{L_{-\frac{1}{2}}}(q,r) = \int_\chi \frac{1}{2}\left[\frac{q(x)+r(x)}{2}+\sqrt{q(x)r(x)}\right]d\mu(x)$$
$$= \frac{1}{2}+\frac{1}{2}D_B(p,q)$$

$$(2.17) \quad D_{L_{-3}}(q,r) = \int_\chi \left[\frac{2r(x)q(x)}{r(x)+q(x)}\cdot r(x)q(x)\right]^{\frac{1}{3}}d\mu(x)$$
$$= \sqrt[3]{2}\int_\chi \frac{r^{\frac{2}{3}}(x)q^{\frac{2}{3}}(x)}{[r(x)+q(x)]^{\frac{1}{3}}}d\mu(x) =: D_{[HG^2]^{\frac{1}{3}}}(p,q)$$

for all $q,r \in \Omega$.

Using Lemma 1, we can state the following fundamental theorem regarding the position of the $L_p$–divergence measures.

**Theorem 1.** *For any $q, r \in \Omega$, we have the inequality*

$$(2.18) \quad 1-\frac{1}{2}D_v(r,q) \leq D_{L_u}(r,q) \leq D_{L_s}(r,q) \leq 1+\frac{1}{2}D_v(r,q)$$

*for all $-\infty \leq u < s \leq \infty$.*
*In particular, we have*

$$(2.19) \quad 1-\frac{1}{2}D_v(p,q) \leq D_{[HG^2]^{\frac{1}{3}}}(p,q) \leq D_B(p,q) \leq D_L(p,q)$$
$$\leq \frac{1}{2}+\frac{1}{2}D_B(p,q) \leq D_I(p,q) \leq 1,$$

*where*

$$D_L(r,q) := \int_\chi \left[\frac{r(x)-q(x)}{\ln r(x)-\ln q(x)}\right]d\mu(x) \text{ is the Logarithmic divergence}$$

*and*

$$D_I(r,q) = \frac{1}{e}\int_\chi \left[\frac{[r(x)]^{r(x)}}{[q(x)]^{q(x)}}\right]^{\frac{1}{r(x)-q(x)}} d\mu(x) \text{ is the Identric divergence.}$$

**Remark 1.** *From (2.18), we can conclude the following inequality for the $L_p$–divergence measure in terms of the variation distance*

$$(2.20) \quad |D_{L_s}(r,q)-1| \leq \frac{1}{2}D_v(r,q),\ r,q \in \Omega$$

*for all $s \in [-\infty,\infty]$. The constant $\frac{1}{2}$ is sharp.*

Indeed, if we assume that (2.20) holds with another constant $c > 0$, i.e.,

$$(2.21) \quad |D_{L_s}(r,q)-1| \leq cD_v(r,q),$$

then, choosing $s = \infty$, we obtain

$$\left|1+\frac{1}{2}D_v(r,q)-1\right| \leq cD_v(r,q) \text{ for all } r,q \in \Omega,$$

which implies that $c \geq \frac{1}{2}$.

For $r \in \mathbb{R}$, we define the $\alpha-th$ *power mean* of the positive numbers $a, b$ by (see [1, p. 133])

$$(2.22) \qquad M^{[\alpha]}(a,b) := \begin{cases} \left(\frac{a^\alpha + b^\alpha}{2}\right)^{\frac{1}{\alpha}} & \text{if } \alpha \neq 0,\ \alpha \neq \pm\infty; \\ \sqrt{ab} & \text{if } \alpha = 0; \\ \max\{a,b\} & \text{if } \alpha = +\infty; \\ \min\{a,b\} & \text{if } \alpha = -\infty. \end{cases}$$

We define the $\alpha-$*power divergence measure* by
$$(2.23)$$
$$D_{M^{[\alpha]}}(p,q) := \begin{cases} \int_\chi \left[\frac{p^\alpha(x)+q^\alpha(x)}{2}\right]^{\frac{1}{\alpha}} d\mu(x) & \text{if } \alpha \neq 0,\ \alpha \neq \pm\infty; \\ \int_\chi \sqrt{p(x)q(x)} d\mu(x) & \text{if } \alpha = 0; \text{ (Bhattacharyya distance)} \\ 1 + \frac{1}{2} D_v(p,q) & \text{if } \alpha = +\infty; \\ 1 - \frac{1}{2} D_v(p,q) & \text{if } \alpha = -\infty. \end{cases}$$

Since $M^{[\alpha]}(a,b) = M^{[\alpha]}(b,a)$ for all $a, b > 0$ and $\alpha \in [-\infty, \infty]$, we can conclude that the $\alpha-$*power divergences are symmetrical*. Now, if we consider the continuous mappings (which are not necessarily convex)

$$(2.24) \qquad f_\alpha(x) := \begin{cases} \left[\frac{x^\alpha + 1}{2}\right]^{\frac{1}{\alpha}} & \text{if } \alpha \neq 0,\ \alpha \neq \pm\infty; \\ \sqrt{x} & \text{if } \alpha = 0,\ x \in (0, \infty); \\ 1 + \frac{1}{2}|x-1| & \text{if } \alpha = +\infty; \\ 1 - \frac{1}{2}|x-1| & \text{if } \alpha = -\infty \end{cases}$$

and taking into account that $M^{[\alpha]}(a,b) = aM^{[\alpha]}\left(1, \frac{b}{a}\right)$, we deduce that

$$\begin{aligned}(2.25) \qquad D_{f_\alpha}(p,q) &= \int_\alpha p(x) f_\alpha\left(\frac{q(x)}{p(x)}\right) d\mu(x) \\ &= \int_\chi p(x) M^{[\alpha]}\left(\frac{q(x)}{p(x)}, 1\right) d\mu(x) \\ &= \int_\chi M^{[\alpha]}(p(x), q(x)) d\mu(x) = D_{M^{[\alpha]}}(p,q)\end{aligned}$$

for all $p, r \in \Omega$, which shows that the $\alpha-$*power divergence measures* can be interpreted as $f-$*Csiszár divergences* for $f = f_\alpha$, which are not necessarily convex.

The following result concerning the fundamental property of the $\alpha-$power means holds (see [1, p. 133 and p. 159]).

**Lemma 2.** *Let $a, b > 0$. Then*

$$(2.26) \qquad M^{[-\infty]}(a,b) \leq M^{[\alpha]}(a,b) \leq M^{[\beta]}(a,b) \leq M^{[+\infty]}(a,b)$$

for $-\infty \leq \alpha < \beta \leq \infty$.
Also,

(2.27) $$\lim_{\alpha \to 0} M^{[\alpha]}(a,b) = G(a,b), \quad \lim_{\alpha \to \pm\infty} M^{[\alpha]}(a,b) = M^{[\pm\infty]}(a,b).$$

In particular

(2.28) $$M^{[-1]}(a,b) = H(a,b).$$

Using the above lemma, we can state the following theorem concerning the location of the $\alpha$–power divergence measure.

**Theorem 2.** *For any $p, q \in \Omega$, we have:*

(2.29) $$1 - \frac{1}{2}D_v(p,q) \leq D_{M^{[\alpha]}}(p,q) \leq D_{M^{[\beta]}}(p,q) \leq 1 + \frac{1}{2}D_v(p,q)$$

*for $-\infty \leq \alpha < \beta \leq \infty$.*
*In particular, we have*

(2.30) $$1 - \frac{1}{2}D_v(p,q) \leq D_{Ha}(p,q) \leq D_B(p,q) \leq 1 + \frac{1}{2}D_v(p,q),$$

*where $D_{Ha}(p,q)$ is the Harmonic divergence and $D_B(p,q)$ is the Bhattacharyya distance.*

**Remark 2.** *From (2.29), we may conclude the following inequalities for the $\alpha$-power divergence measures in terms of the variational distance*

(2.31) $$|D_{M^{[\alpha]}}(p,q) - 1| \leq \frac{1}{2}D_v(p,q)$$

*for any $p,q \in \Omega$ and $\alpha \in [-\infty, \infty]$ and the constant $\frac{1}{2}$ is sharp.*

In what follows, by the use of a result by Pittenger [1, p. 349], we obtain inequalities that are related to logarithmic and power means:

**Theorem 3.** *Let $a, b > 0$, $-\infty < r < \infty$ and define*

$$\begin{aligned}
r_1 &= \min\left\{\frac{r+2}{3}, r \cdot \frac{\ln 2}{\ln r + 1}\right\}, r > -1, r \neq 0 \\
&= \min\left\{\frac{2}{3}, \ln 2\right\}, r = 0 \\
&= \min\left\{\frac{r+2}{3}, 0\right\}, r \leq -1,
\end{aligned}$$

*with $r_2$ as defined above, but with $\max$ instead of $\min$, then*

(2.32) $$M^{[r_1]}(a,b) \leq L_r(a,b) \leq M^{[r_2]}(a,b),$$

*with equality iff $a = b$ or $r = 1, -\frac{1}{2}$ or $-2$. The values $r_1$ and $r_2$ are sharp.*

We are able to establish the following relationship between the power divergence and the generalized logarithmic divergence.

**Theorem 4.** *For any $p, q \in \Omega$, we have:*

(2.33) $$D_{M^{[r_1]}}(p,q) \leq D_{L_r}(p,q) \leq D_{M^{[r_2]}}(p,q),$$

*where $r_1, r_2$ are as defined above.*

**Remark 3.** a) *In the case $r = -1$, the corresponding inequality in (2.32) is due to Lin (see [1, p. 349]) and for this $r$ we obtain from (2.33)*

$$(2.34) \qquad B(p,q) \leq D_L(p,q) \leq D_{M^{\left[\frac{1}{3}\right]}}(p,q), \quad p,q \in \Omega.$$

b) *If $-2 < r < -\frac{1}{2}$ or $r > 1$, then $r_2 = \frac{r+2}{3}$ and the right hand side of (2.33) becomes*

$$(2.35) \qquad D_{L_r}(p,q) \leq D_{M^{\left[\frac{r+2}{3}\right]}}(p,q), \quad p,q \in \Omega.$$

Using the above means, we can imagine all the sets of other divergences that can be constructed by the use of different contributions of these means. For example, we can define for $p, q \in \Omega$

$$D_{(AG)^{\frac{1}{2}}}(p,q) := \int_\chi \sqrt{A(p(x),q(x))\, G(p(x),q(x))}\, d\mu(x);$$

$$D_{(LI)^{\frac{1}{2}}}(p,q) := \int_\chi \sqrt{L(p(x),q(x))\, I(p(x),q(x))}\, d\mu(x)$$

or even

$$D_{(GI)^{\frac{1}{2}}}(p,q) := \int_\chi \sqrt{G(p(x),q(x))\, I(p(x),q(x))}\, d\mu(x).$$

Using Alzer's result for means (see [1, p. 350])

**Theorem 5.** *If $a, b > 0$, we have*

$$\begin{aligned}
\sqrt{A(a,b)\, G(a,b)} &< \sqrt{L(a,b)\, I(a,b)} < M^{\left[\frac{1}{2}\right]}(a,b); \\
L(a,b) + I(a,b) &< A(a,b) + G(a,b); \\
\sqrt{G(a,b)\, I(a,b)} &< L(a,b) < \frac{1}{2}\left[G(a,b) + I(a,b)\right],
\end{aligned}$$

and we may state the following theorem concerning the above divergence measures.

**Theorem 6.** *For any $p, q \in \Omega$, we have*

$$(2.36) \qquad D_{(AG)^{\frac{1}{2}}}(p,q) < D_{(LI)^{\frac{1}{2}}}(p,q) < D_{M^{\left[\frac{1}{2}\right]}}(p,q),$$

$$(2.37) \qquad D_L(p,q) + D_I(p,q) < 1 + B(p,q),$$

$$(2.38) \qquad D_{(GI)^{\frac{1}{2}}}(p,q) < D_I(p,q) < \frac{1}{2}\left[B(p,q) + D_I(p,q)\right].$$

**Remark 4.** *In this way, we have shown that **any result for special means** can be imported for the divergence measure generated by these means, providing a very rich universe of facts in comparing the new distances with the other distances which have already become classics in Information Theory: Bhattacharyya distance $B(p,q)$, Harmonic distance $Ha(p,q)$ or variation distance $D_v(p,q)$, etc.*

School of Communications and Informatics, Victoria University of Technology, PO 14428, Melbourne City MC, Victoria, 8001, Australia
*E-mail address*: sever.dragomir@vu.edu.au